\documentclass{amsproc}

\usepackage{a4wide}
\usepackage{amsmath}
\usepackage{enumerate}

\usepackage{citehack}
\usepackage{amsmath,amsthm}
\usepackage{amsfonts}
\usepackage{amssymb}
\usepackage{longtable}
\usepackage[matrix,arrow,curve]{xy}

\usepackage{lipsum}

\newtheorem{theorem}{Theorem}
\newtheorem{cor}{Corollary}
\newtheorem{prop}{Proposition}

\newtheorem{ex}{Example}{\bf}{\rm}

\def\Real{\mathbb{R}}

\def\U{\mathcal{U}}
\def\V{\mathcal{V}}

\def\D{\mathcal{D}}
\def\C{\mathcal{C}}
\def\M{\mathcal{M}}

\def\F{\mathcal{F}}
\def\vH{{\check{H}}}

\def\GLn{\text{\rm GL}(n,\Real)}
\def\On{\text{\rm O}(n)}
\def\GL1{\text{\rm GL}(1,\Real)}
\def\O1{\text{\rm O}(1)}

\def\id{\mathop\text{\rm id}\nolimits}
\def\pt{\mathop\text{\rm pt}\nolimits}

\def\WOn{\underline{\text{\rm WO}}_n}
\def\Wn{\underline{\text{\rm W}}_n}
\def\WGLn{\underline{\text{\rm WGL}}_n}

\usepackage{stackrel}

\newcommand{\be}{\begin{equation}}
\newcommand{\ee}{\end{equation}}

\begin{document}


\title{Comparison  of  approaches to characteristic classes of foliations}

\author{Anton S. Galaev}\thanks{$^1$University of Hradec Kr\'alov\'e, Faculty of Science, Rokitansk\'eho 62, 500~03 Hradec Kr\'alov\'e,  Czech
	Republic\\
	E-mail: anton.galaev(at)uhk.cz}

\maketitle

\begin{abstract} It is shown that the characteristic classes of foliations that were defined by Losik and that take values in the de~Rham
	cohomology of the space of infinite order frames  over the leaf
	space may be mapped to the characteristic classes with values in 
	the \v{C}ech-de~Rham cohomology of the leaf space studied in details by
	Crainic and Moerdijk. This map is in general non-injective. All constructions are done using Losik's approach to Gelfand formal geometry. A similar result is obtained for the exotic characteristic classes as well as for the group actions of the diffeomorphisms. As illustrating examples, foliations of codimension one are discussed. 



\end{abstract}

\section*{Introduction}

The theory of characteristic classis of foliations starts with the
discovery of the Godbillon-Vey class \cite{GV} for a foliation of
codimension one. Berstein and Rozenfeld \cite{BR}, and Bott and Haefliger \cite{BH}
defined characteristic classes for foliations of arbitrary
codimension $n$. These classes are given by the generators of the
relative Gelfand-Fuchs cohomology $H^*(W_n,\On)$ of the Lie algebra
of formal vector fields on $\Real^n$ and take values in the
cohomology $H^*(M)$ of the underlying manifold. The map
$H^*(W_n,\On)\to H^*(M)$ may be factored through the cohomology
$H^*(B\Gamma_n)$ of the classifying space $B\Gamma_n$ of the
Haefliger groupoid $\Gamma_n$.

Losik~\cite{L90,L94,L15} introduced the notion of a generalized
atlas on the leaf space $M/\F$ of a foliation $\F$ on a smooth
manifold $M$ showing that there is a reach smooth structure on
$M/\F$. Spaces with such structures are called $\D_n$-spaces.  In
particular, this approach allowed him to define the characteristic
classes as elements of the cohomology $H^*(S(M/\F)/\On)$, where
$S(M/\F)$ is the bundles of frames of infinite order over $M/\F$.
These classes may be projected to the usual ones in $H^*(M)$.

Following Haefliger, one may construct the characteristic map $H^*(W_n,\On)\to H^*(B\Gamma_n) \to H^*(BG)$ and obtain characteristic classes with values in the cohomology of the classifying space of the holonomy groupoid of the foliation \cite{Bott,BSS,Hae}. Again,
these classes may be projected to the usual ones in $H^*(M)$. In the more recent work \cite{CM}, Crainic and Moerdijk   studied in details the
characteristic classes of foliations as elements of the
\v{C}ech-de~Rham cohomology $\vH^*(M/\F)$ of the leaf space $M/\F$. The leaf space $M/\F$ was considered in terms of the holonomy groupoid.   There is an isomorphism  
$\vH^*(M/\F)\cong H^*(BG)$ \cite{CM}.

The aim of the present paper is to compare the approaches from the last two paragraphs.
Using Losik's theory, we recover the construction of  Crainic and
Moerdijk, and we  show that Losik's characteristic classes may be
mapped to these of Crainic and Moerdijk. Let us sketch the idea. To
a $\D_n$-space $X$ is associated an \'etale groupoid $G_X$ such that
there exists an isomorphism $\vH^*(X)\cong H^*(BG_X)$ \cite{L15}.
The category of $\D_n$-spaces possesses the terminal object, which
is a point $\pt$ with an additional structure. The associated
groupoid to that space is exactly the Haefliger groupoid $\Gamma_n$.
For a $\D_n$-space, the unique map $X\to\pt$ induces the
homomorphism $H^*(B\Gamma_n)\to \vH^*(X)$, which in the case
$X=M/\F$ was defined by Crainic and Moerdijk. Next, there is a map
from the de~Rham cohomology $H^*(S(X)/\On)$ to $\vH^*(S(X)/\On)$; it
holds $H^*(S(\pt)/\On)\cong H^*(W_n,\On)$; and we prove that
$\vH^*(S(X))\cong \vH^*(X)$. This gives us the sequence of the maps
$$H^*(W_n,\On)\cong H^*(S(\pt)/\On)\to H^*(S(X)/\On)\to
\vH^*(S(X)/\On)\cong\vH^*(X),$$ which provides the characteristic
homomorphism of Crainic and Moerdijk for the case $X=M/\F$. In
particular, we see that the Losik classes living in $H^*(S(X)/\On)$
are mapped to the classes contained in $\vH^*(M/\F)$.
 If we take $X=\pt$, then we obtain the new construction of the Bott map
$H^*(W_n,\On)\to H^*(B\Gamma_n)$.

Let $n=1$. The Godbillon-Vey class is defined by a generator of
$H^3(W_1,\O1)\cong\Real$. Let $\F$ be the Reeb foliation on the
three dimensional sphere. It is known that the usual Godbillon-Vey
class of $\F$ is trivial in $H^3(S^3)$ \cite{Tamura}. The triviality of the Godbillon-Vey class in  $H^3(S(S^3/\F)/\O1)\cong H^3(BG)$ may be shown by means of the noncommutative geometry: if the Godbillon-Vey class of a foliation is non-trivial in $H^3(BG)$, then the corresponding von~Neumann algebra has a non-trivial type III component \cite{H-k}, \cite[p. 245, p. 261]{Connes}. On the other hand, the von~Neumann algebra of the Reeb foliation is of type I, \cite[p. 54]{Connes}. In contrast, in \cite{BGG}, we show
that the Godbillon-Vey class for some Reeb foliations $\F$ is non-trivial in
$H^3(S(S^3/\F)/\O1)$. We conclude that the map $H^*(S(M/\F)/\On)\to \vH^*(M/\F)$ that we constructed here is in general non-injective. 

The above construction gives also the maps $H^*(W_n,\On)\to H^*(S(M/G)/\On)\to H^*(G;\Omega(M))$ for a group $G$ acting by diffeomorphisms on an $n$-dimensional manifold $M$. For $M=S^1$ this allows to reconstruct the Bott-Thurston cocycle, which is an element of $H^*(G;\Real)$. Using results form \cite{BGG}, an example of $G$ acting on $S^1$ is given with the Godbillon-Vey class trivial in    
$H^*(G;\Real)$ and non-trivial in $H^*(S(S^1/G)/\O1)$.

 For foliations $\F$
with the trivializable bundle of frames of the normal bundle one may consider characteristic classes
are given by a map $H^*(W_n)\to H^*(M)$. Losik considered these
classes with values in $H^*(S(M/\F))$. In a way as above,
for a $\D_n$-space $X$ with the trivializable frame bundle   we
construct homomorphisms $$H^*(W_n)\to H^*(S(X))\to \vH^*(X),\quad
H^*(W_n)\to H^*(B\bar\Gamma_n)\to \vH^*(X)$$ defining the so-called exotic
characteristic classes with values in the \v{C}ech-de~Rham
cohomology. Here $\bar\Gamma_n$ is the groupoid with the space of
object the frame of bundles over $\Real^n$ and morphisms the germs
of the extensions of the local diffeomorphisms of $\Real^n$. In
particular, we get a new construction of the Bott map $H^*(W_n)\to
H^*(B\bar\Gamma_n)$.

Losik  considered also the Chern classes of foliations that are
defined by a map $H^*(W_n,\GLn)\to H^*(S(M/\F))$. Compared with the
elements from $H^*(W_n,\On)$, the elements $H^*(W_n,\GLn)$ do not
give any new classes in $H^*(M)$ and $\vH^*(M/\F)$. On the other
hand, elements of $H^*(W_n,\GLn)$ define classes in $H^*_\F(M)$
similar to the Vey class \cite{CCII,H-H,Hurder00}. Let $n=1$. Denote by $c_1$
a generator of $H^2(W_1,\GL1)\cong\Real$. The class defined by this
generator is always trivial in $H^2(M)$ and $\vH^2(M/\F)$. This
generator defines the Vey class in $H^2_\F(M)$. Consider again the
Reeb foliation. The Vey class for that foliation is trivial. In
contrast, Losik \cite{L15} proved (see also \cite{Baz-Gal}) that $c_1$ defines non-trivial class in $H^2(S(S^3/\F)/\GL1)$. We thus may  conclude  that Losik's
classes are most informative. In particular, the first Chern class indicates that the Reeb foliation contains a compact leaf with a non-trivial holonomy \cite{Baz-Gal} (the Godbillon-Vey class is zero whenever only compact leaves have non-trivial holonomy \cite{Ref82}). The triviality of the Godbillon-Vey class with values in   $H^3(S(S^3/\F)/\O1)$ depends on the choice of the functions defining the Reeb foliation \cite{BGG} and it may provide more delicate information about  growth rate of the leaves near the compcat leaf.

In what follows, all cohomology are considered with real
coefficients.

\section{Classical approaches}\label{sec1}

Let $W_n$ denote the Lie
algebra  of formal vector fields on $\Real^n$. First we recall the description of the
Gelfand-Fuchs cohomology of $W_n$ \cite{Fuks}.
 Consider the  differential graded algebra $\Wn$
with the generators $$y_1,\dots,y_n,c_1,\dots,c_n$$ of degrees ${\rm
deg} y_i=2i-1$, ${\rm deg} c_i=2i$, the relations
$$y_iy_j=-y_jy_i,\quad y_ic_j=c_jy_i,\quad c_ic_j=c_jc_i,$$
$$c_{i_1}\cdots c_{i_s}=0\quad\text{if}\quad i_1+\cdots +i_s>n$$ and the differential satisfying $dy_i=c_i$,
$dc_i=0$. The cohomology $H^*(W_n)$ are isomorphic to $H^*(\Wn)$.
The relative cohomology $H^*(W_n,\On)$ are isomorphic to the cohomology
$H^*(\WOn)$
 of the subcomplex generated by the elements
$y_1,y_3,y_5,\dots,c_1,\dots, c_n$; the cohomology $H^*(W_n,\GLn)$
 are isomorphic to the cohomology $H^*(\WGLn)$ of the
subcomplex generated by $c_1,\dots, c_n$.

Now we recall shortly the classical definitions of the
characteristic classes of foliations \cite{BR,Bott,BH,F73}.
Foliations $\F$ of codimension $n$ on  smooth manifolds are
classified by the homotopy classes of  maps \be\label{fF}f_\F:M\to B\Gamma_n,\ee where
$B\Gamma_n$ is the classifying space of the Haefliger groupoid
$\Gamma_n$. In this way one obtains the characteristic map
$$f^*_\F:H^*(B\Gamma_n)\to H^*(M).$$
Using the Chern-Weil theory, one may construct the map
$$k_\F: H^*(\WOn)\to H^*(M).$$
Using the maps $k_\F$, Bott \cite{Bott} constructed the  map
\be\label{Bk} k: H^*(\WOn) \to H^*(B\Gamma_n).\ee  It holds
$\lambda_\F=f^*_\F \circ k$. The images of the generators of
$H^*(\WOn)$ under the map $\lambda_\F$ define {\it characteristic
classes} of the foliation $\F$.

Let $P(\F)$ denote the frame bundle of the normal bundle of $\F$.
Suppose that the bundle $P(\F)$ is trivial, i.e., there exists a
section $s:M\to P(\F)$. Then there is a classifying  map $$\bar
f_\F:M\to B\bar\Gamma_n,$$ where $\bar\Gamma_n$ is the groupoid
whose objects form the frame bundle $P(\Real^n)$ over $\Real^n$, and
the morphisms are the germs of the extensions to $P(\Real^n)$ of the
local diffeomorphisms of $\Real^n$ (the classifying space
$B\bar\Gamma_n$ is the homotopy-theoretic fiber $F\Gamma_n$ of the
map $B\Gamma_n\to B{\rm GL}(n,\Real)$). There is a map
\be\label{Bbk} \bar k: H^*(\Wn) \to H^*(B\bar\Gamma_n)\ee
constructed in \cite{Bott}. The images of the generators of
$H^*(\Wn)$ under the map $$\bar f^*_\F \circ \bar k: H^*(\Wn)\to
H^*(M)$$ define the so-called exotic characteristic classes of
the foliation $\F$.

 Alternatively one may proceed in the following way \cite{BR,F73}.
Let $S(\F)$ be the space of jets of infinite order at zero of
submersions from $M$ to $\Real^n$ that are constant on the leaves of
$\F$. Let $S'(\F)=S(\F)/\GLn$ and $S''(\F)=S(\F)/\On$. On each of
the spaces $S(\F)$, $S'(\F)$, $S''(\F)$, there is a canonical 1-from
(the Gelfand-Kazhdan form) with values in $W_n$ that delivers the
homomorphisms
$$H^*(W_n)\to H^*(S(\F))=H^*(P(\F)),$$
$$H^*(W_n,\GLn)\to
 H^*(S'(\F))=H^*(M),$$
$$H^*(W_n,\On)\to H^*(S''(\F))=H^*(M).$$ If the bundle
$P(\F)$ is trivializable, then one gets the map
$$H^*(W_n)\to H^*(P(\F))\to H^*(M).$$ This gives another
construction of the  characteristic classes.
The map
$$H^*(W_n,\GLn)\to
 H^*(M)$$ does not give any new classes, since there is the
 following commutative diagram:
$$\xymatrix{
H^*(W_n,{\rm GL}(n,\Real)) \ar@{->}[rr]  \ar@{->}[d]
& & \ar@{->}[d]H^*(S'(\F))\\
H^*(W_n,{\rm O}(n)) \ar@{->}[rr]    & & H^*(S''(\F)) }$$
 Below we will see that the generators of $H^*(W_n,\GLn)$ define
 certain interesting classes with values in the cohomology of a bundle
 over the leaf space of the foliation.

Finally note that  the Godbillon-Vey class of a codimension one
foliation  is defined by the class $[y_1c_1]\in H^3(\underline{\rm
WO}_1)=H^3(W_1,\rm{O}(1)).$

\section{Characteristic classes with values in $\vH^*(M/\F)$}\label{sec2}

Let us briefly recall the construction of the  \v{C}ech-de~Rham
cohomology following Crainic and Moerdijk~\cite{CM}. Consider a foliation
$\F$ of codimension $n$ on a smooth manifold $M$. Let $\U$ be a
family of transversal sections of $\F$. Such a family is called a
transversal basis if for each transversal section $V$ of $\F$ and
each point $y\in V$, there exists a section $U\in\U$ and a holonomy
embedding $h:U\to V$ such that $y\in h(U)$. Consider the double
complex \be\label{Cpq}
C^{p,q}=\prod_{U_0\stackrel{h_1}{\longrightarrow}\cdots
\stackrel{h_p}{\longrightarrow} U_p}\Omega^q(U_0),\ee where the
product ranges over all $p$-tuples of holonomy embeddings between
transversal sections from a fixed transversal basis $\U$. The
vertical differential is defined as $$(-1)^p d:C^{p,q}\to
C^{p,q+1},$$ where $d$ is the usual de~Rham differential. The
horizontal differential $$\delta:C^{p,q}\to C^{p+1,q}$$ is given by
\begin{multline}\label{delta}
(\delta\omega)(h_1,\dots,h_{p+1})=h_1^*\omega(h_2,\dots,h_{p+1})+\\
\sum_{i=1}^p(-1)^i\omega(h_1,\dots,h_{i+1}h_i,\dots,h_{p+1})+(-1)^{p+1}\omega(h_1,\dots,
h_p).
\end{multline}
The cohomology of this complex is called the \v{C}ech-de~Rham
cohomology of the leaf space $M/\F$ with respect to the cover $\U$
and is denoted by
$$\vH^*_\U(M/\F).$$
A complete transversal basis $\U$ may be obtained from a foliation
atlas $\tilde \U$ of $M$. This defines a map of the complexes
$C^{p,q}(\U)\to C^{p,q}(\tilde \U)$, which induces the map
\be\label{prH}\vH^*_\U(M/\F)\to \vH^*_{\tilde\U}(M)\cong H^*(M),\ee
where $H^*(M)$ is the de~Rham cohomology of $M$. Next, there are
natural isomorphisms \be\label{isom1} \vH^*_\U(M/\F)\cong H^*(B{\rm
Hol}(M,\F))\cong H^*(B{\rm Hol}_T(M,\F)),\ee where $T$ is a complete
transversal, $\U$ is a basis of $T$, $B{\rm Hol}(M,\F)$ is the
classifying space of the holonomy groupoid, and $B{\rm Hol}_T(M,\F)$
is the classifying space of the holonomy groupoid restricted to $T$.

Using the Chern-Weil theory, Crainic and Moerdijk constructed the
characteristic map \be\label{vkF} \check k_\F:
H^*(\WOn)\to\vH^*_\U(M/\F). \ee Even more, for a given \'etale
groupoid $G$, they defined the \v{C}ech-de~Rham cohomology
$\vH_\U(G)$ and constructed the map \be\label{vkG} \check k_G:
H^*(\WOn)\to\vH^*_\U(G)\cong H^*(BG). \ee For the case $G=\Gamma_n$
this gives the  new construction of the map \eqref{Bk}. It is shown
that the classifying map \eqref{fF} induces the characteristic map
\be\label{CM19}  H^*(B\Gamma_n)\to\vH^*_\U(M/\F). \ee
 Thus
there is the following commutative diagram:
$$\xymatrix{ H^*(\WOn) \ar@{->}^k[rr]
\ar@{->}_{\check k_\F}[rd]& &H^*(B\Gamma_n)\ar@{->}[ld]\\
& \vH^*_\U(M/\F) &}$$ The images of the generators of $H^*(\WOn)$
under the map $\check k_\F$ are characteristic classes of the
foliation $\F$ living in $\vH^*_\U(M/\F)$. The map \eqref{prH} sends
these classes to the characteristic classes considered in
Section~\ref{sec1}.

\section{Losik's approach}\label{sec3}

The main idea of Losik's approach to the leaf spaces of foliations is
to introduce a notion of a generalized smooth atlas on such spaces
in order to apply to them the technics  from the theory of smooth
manifolds \cite{L90,L94,L15}. In this approach, the main notion is a
$\D_n$-space, where $\D_n$ is the category whose objects are open
subsets of $\Real^n$, and morphisms are \'etale (i.e., regular)
maps. The dimension $n$ may be infinite; in this case we use the
definitions from \cite{BR} of manifolds with the model space
$\Real^\infty$.

Let us recall the definition of a $\D_n$-space. Let $X$ be a set. A
$\D_n$-chart on $X$ is a pair $(U,k)$, where $U\subset \Real^n$ is
an open subset, and $k:U\to X$ is an arbitrary map. For two charts
$k_i:U_i\to X$, a morphism of charts is an \'etale map $m:U_1\to
U_2$ such that $k_2\circ m=k_1$. Let $\Phi$ be a set of charts and
let $\C_\Phi$ be a category whose objects are elements of $\Phi$ and
morphisms are some morphisms of the charts. The set $\Phi$ is called
a $\D_n$-atlas on $X$ if $X=\varinjlim J$, where $J:\C_\Phi\to{\rm
Sets}$ is the obvious functor. A $\D_n$-space is a set $X$ with a
maximal $\D_n$-atlas. An atlas $\Phi$ is called full if for each
chart $(V,l)$  from the corresponding maximal atlas $\bar \Phi$ and
each point $y\in V$, there exists a chart $(U,k)\in\Phi$ and a
morphism $m:(U,k)\to (V,l)$ such that $y\in m(U)$.

A map $f:X\to Y$ of $\D_n$-spaces is called a morphism of
$\D_n$-spaces if, for any  $\D_n$-chart $k$ from the atlas on $X$,
$f \circ k$ is a $\D_n$-chart from the maximal atlas on $Y$.

If $\F$ is a foliation of codimension $n$ on a smooth manifold $M$,
then the leaf space $M/\F$ is a $\D_n$-space. The maximal
$\D_n$-atlas on $M/\F$ consists of the projections $U\to M/\F$,
where $U$ is a transversal which is the embedded to $M$ open subset
of $\Real^n$. These transversals may be obtained from a foliation
atlas on $M$. A full atlas may be obtained from a complete
transversal as in the previous section.

A full atlas $\Phi$ of a $\D_n$-space $X$  gives rise to a smooth
groupoid $ G_X$. The set of objects of $G_X$ is the union of the
domains of the charts form $\Phi$, and the morphisms are germs of
the morphisms of charts. The groupoid $G_X$ may be reduced
essentially, if there is a surjection $(G_X)_0\to M$ to a smooth
manifold $M$. For the reduced groupoid $G_X$, $(G_X)_0= M$, and the
elements of $(G_X)_1$ are the germs of local diffeomorphisms of $M$
that can be lifted to the morphisms form $\C_\Phi$. If $X=M/F$, then
the reduced groupoid $G_X$ coincides with the holonomy groupoid.

Generally $\D_n$-spaces are orbit spaces of  pseudogroups of local
diffeomorphisms of smooth manifolds. Considering the space $\Real^n$
and the pseudogroup of all local diffeomorphisms of open subsets of
$\Real^n$, we see that the point $\pt$ is a $\D_n$-space. The atlas
of $\pt$ consist of all pairs $(U,k)$, where $U\subset\Real^n$ is an
open subset and $k:U\to \pt$ is the unique map. It is important to
note that $\pt$ is the terminal objects in the category of
$\D_n$-spaces. The reduced groupoid corresponding to the
$\D_n$-space $\pt$ is exactly the Haefliger groupoid $\Gamma_n$.

Each (co)functor from the category   $\D_n$ to the category of sets
may be extended to a (co)functor form the category of $\D_n$-spaces.
In this way one obtains, e.g., the de~Rham complex $\Omega^*(X)$ of
a $\D_n$-space $X$, which defines the de~Rham cohomology $H^*(X)$ of
$X$. E.g., if $X=M/\F$, then $\Omega^*(X)$ coincides with the
complex of basic forms.

One may consider also the category $\M_n$ of $n$-dimensional
manifolds with the \'etale maps as morphisms. Then one gets the
notion of an $\M_n$-space. The categories of $\D_n$-spaces and
$\M_n$-spaces are equivalent. We do not restrict the attention to
$\D_n$-spaces, since now we will consider $\M_\infty$-atlases on the
spaces of frame of infinite order.

Consider the functor $S$ assigning to each open subset
$U\subset\Real^n$ the space of frames of infinite order, i.e., the
space of jets at $0\in\Real^n$ of regular maps from $\Real^n$ to
$U$. Then for each $\D_n$-space $X$, we obtain the space of frames
of infinite order $S(X)$. This space is a $\M_\infty$-space. Each
chart $U\to X$ defines the char $S(U)\to S(X)$, and each morphism of
charts $h:U\to V$ defines the morphism $ S(h):S(U)\to S(V)$. In this
way we get an $\M_\infty$-atlas on $S(X)$. Similarly, let
$S'(U)=S(U)/\GLn$ and $S''(U)=S(U)/\On$. We obtain the
$\M_\infty$-spaces $S'(X)$ and $S''(X)$ with the atlases similar to
the above one.

Let us consider the point $\pt$ as a $\D_n$-space. Each of the
spaces $S(\pt)$, $S'(\pt)$, and $S''(\pt)$ consists of a single
point, on the other hand, the complexes $\Omega^*(S(\pt))$,
$\Omega^*(S'(\pt))$, and $\Omega^*(S''(\pt))$, are naturally
isomorphic to the complexes $C^*(W_n)$, $C^*(W_n,\GLn)$, and
$C^*(W_n,\On)$, respectively.

Let now $X$ be a $\D_n$-space. The unique morphism
$$p_X:X\to\pt$$ of $\D_n$-spaces induces the morphisms
$$S(X)\to S(\pt),\quad S'(X)\to S'(\pt),\quad S''(X)\to S''(\pt)$$
and the characteristic morphisms
$$\chi:H^*(W_n)\cong H^*(S(pt))\to H^*(S(X)),$$ $$\chi':H^*(W_n,\GLn)\cong H^*(S'(pt))\to
H^*(S'(X)),$$ $$\chi'':H^*(W_n,\On)\cong H^*(S''(pt))\to
H^*(S''(X)).$$ The images of the generators of the cohomology under
these maps give Losik's characteristic classes. In the case
$X=M/\F$, these classes may be projected to the characteristic
classes from Section~\ref{sec1} as  follows. The projection $p:M\to
M/\F$ defines the map $S(p):S(\F)\to S(M/\F)$ and two similar maps
$S'(p)$ and $S''(p)$, and one obtains the following chains of
homomorphisms: \begin{align}\label{map1}
H^*(W_n)&\stackrel{\chi}\longrightarrow
H^*(S(M/\F))\stackrel{S(p)^*}\longrightarrow
H^*(S(\F))=H^*(P(\F)),\\ \label{map2}
H^*(W_n,\GLn)&\stackrel{\chi'}\longrightarrow
H^*(S'(M/\F))\stackrel{S'(p)^*}\longrightarrow H^*(S'(\F))=H^*(M),\\
\label{map3} H^*(W_n,\On)&\stackrel{\chi''}\longrightarrow
H^*(S''(M/\F))\stackrel{S''(p)^*}\longrightarrow
H^*(S''(\F))=H^*(M).\end{align}

\section{Comparison of the approaches}\label{sec4}

Generalizing the construction from  \cite{CM}, Losik \cite{L15}
defined the \v{C}ech-de~Rham cohomology for an $\M_n$-space $X$. Let
$\Phi$ be a full $\M_n$-atlas on $X$. The disjoint union of the
domains of the charts from $\Phi$ is an $n$-dimensional manifold.
Let $\U$ be a base of the topology on this manifold. We will refer
to $\U$ as to a complete cover of $X$. The word ,,complete''
stresses that the cover is obtained from a full atlas. Consider the
double complex \be\label{Cpq2}
C_\U^{p,q}(X)=\prod_{U_0\stackrel{h_1}{\longrightarrow}\cdots
\stackrel{h_p}{\longrightarrow} U_p}\Omega^q(U_0),\ee where product
is taken over the strings of composable arrows from $\C_\Phi$. The
differentials of this double complex are defined in the same way as
for the double complex \eqref{Cpq} above. The obtained cohomology are
called the \v{C}ech-de~Rham cohomology of the $\M_n$-space defined
by the cover $\U$ and are denoted by $$\vH^*_\U(X).$$ If $X=M/\F$ is
a leaf space and $\U$ is as in Section \ref{sec1}, then both
complexes coincide.

Losik proved\footnote{In the proof it is implicitly assumed that the
cover consists of contactable sets. To avoid this assumption one may
use the proof of the isomorphism \eqref{isom1} from \cite{CM}. For
that it is enough to apply Theorem~3 from \cite{CM} to the \'etale
groupoid $G_X$.} the existence of the natural  isomorphism
\be\label{lab*}\vH^*_\U(X)\cong H^*(BG_X),\ee where $B G_X$ is the
classifying space of the groupoid $G_X$ defined above. This implies
that the cohomology $\vH^*_\U(X)$ does not depend on the choice
neither of a full $\D_n$-atlas $\Phi$, nor on the base $\U$.

 In
particular, for the $\D_n$-space $\pt$ we get \be\label{vHpt}
\vH^*_\V(\pt)\cong H^*(B\Gamma_n),\ee where $\V$ is a base of the
topology on $\Real^n$.

Now we are going to describe the relation between  the
characteristic classes defined by Crainic and Moerdijk with the
classes of Losik.

Let $X$ be a $\D_n$-space with a full atlas $\Phi$. A cover $\U$ for
$X$ defines the cover $\U_S$ for $S(X)$ consisting of the elements
$S(U)$. Each morphism of charts $h:U\to V$ defines the morphism
$S(h): S(U)\to S(V)$. We consider in the definition of the
\v{C}ech-de~Rham cohomology for $S(X)$ only such morphisms. We use
similar notations for the spaces $S'(X)$, $S''(X)$ and $P(X)$, where
$P(X)$ is the frame bundle of $X$.

\begin{prop}\label{prop1}  The projection $\pi: S''(X)\to X$ induces the
isomorphism $$\pi^*:  \vH^*_\U(X)\to \vH^*_{ \U_{S''}}(S''(X)).$$
\end{prop}

{\bf Proof.} The map of the complexes $$\pi^*: C^{p,q}_\U(X)\to
C^{p,q}_{\U_{S''}}(S''(X))$$ is defined by the equality
$$(\pi^*\omega)( S''(h_1),\dots, S''(h_p))=\pi^*(\omega(h_1,\dots,h_p)),$$
where $\pi^*$ on the right hand side is the map acting on the
differential forms, and the chain of arrows
$$U_0\stackrel{h_1}{\longrightarrow}\cdots
\stackrel{h_p}{\longrightarrow} U_p$$ is the one that uniquely
defines the chain
$$ S''(U_0)\stackrel{ S''(h_1)}{\longrightarrow}\cdots
\stackrel{ S''(h_p)}{\longrightarrow}  S''(U_p).$$ It is clear that
the map $\pi^*$ induces the isomorphisms
$$\pi^*:H^q(C^{p,\bullet}_\U(X))\to H^q(C^{p,\bullet}_{\U_{S''}}(S''(X))).$$
This and Theorem 1.1 from \cite{DHN} imply the proof of the proposition.
\qed

{\bf Remark.} It is known that for a manifold $M$, the spaces
$S''(M)$ and $M$ are weakly homotopy equivalent and consequently
have isomorphic de~Rham cohomology. This is not the case for a
$\D_n$-space $X$, consider, e.g., $X=\pt$.

\smallskip

For a $\D_n$-space $X$, one defines the natural map (cf. \cite[Sec.
6.2]{CM})
$$j:\Omega^*(X)\to C^{0,*}_\U(X)$$ as follows:
$$j(\omega)(U_0)=\omega|_{U_0},$$ here $U_0$ corresponds in
\eqref{Cpq2} to the strings of length zero. This defines the
homomorphism of cohomology \be\label{morj}
j:H^*(X)\to\vH^{*}_\U(X).\ee Taking, $X=\pt$, we immediately obtain
the new construction of the map \eqref{Bk}:
\be\label{newk}k:H^*(\WOn)\cong H^*(W_n,\On)\cong
H^*(S''(\pt))\stackrel{j}{\longrightarrow}\vH^*_{\V_{S''}}(S''(\pt))\cong\vH^*_\V(\pt)\cong
H^*(B\Gamma_n).\ee

Let $X$ be a $\D_n$-space with an atlas $\Phi$. It is clear that the
domain of each chart from $\Phi$ is a domain of the chart on $\pt$,
and any morphism of charts from $\Phi$ is a morphism of the
corresponding charts on $\pt$. This shows that the unique map
$$p_X:X\to\pt$$ induces a homomorphism of cohomology
$$p^*_X:\vH_\V^*(\pt)\to \vH^*_\U(X).$$
Together with the isomorphism \eqref{vHpt} this gives the map
\be\label{charmap} H^*(B\Gamma_n)\cong\vH_\V^*(\pt)\to
\vH^*_\U(X).\ee Using the just constructed map $k$, we instantly
obtain the characteristic homomorphism \be\label{newkX}
\check{k}_X:H^*(\WOn)\to \vH^*_\U(X).\ee

\begin{theorem}\label{Thcoincide}
The homomorphism \eqref{newk} coincides with the homomorphism
\eqref{Bk}. If $X=M/\F$ is the leaf space of a foliation, then the
homomorphism \eqref{charmap} coincides with the homomorphism
\eqref{CM19}; consequently, the homomorphism \eqref{newkX} coincides
with the homomorphism \eqref{vkF}.\end{theorem}

{\bf Proof.} Suppose that $X=M/\F$ is the leaf space of a foliation.
A morphism of $\D_n$ spaces $f:X\to Y$ induces the morphism of the
corresponding groupoids $\hat f:G_X\to G_Y$  and the morphism $\hat
f:BG_X\to BG_Y$ of the classifying spaces \cite{L15}. By the
naturality of the isomorphism \eqref{lab*}, the map $p_X:X\to\pt$
delivers the commutative diagram
$$\xymatrix{
\vH^*_\V(\pt) \ar@2{-}[r]
\ar@{->}^{p_X^*}[d]& H^*(B\Gamma_n)\ar@{->}^{\hat p^*_X}[d]\\
\vH^*_\U(X) \ar@2{-}[r] & H^*(BG_X)}$$ Since $X=M/F$, $G_X$ is the
holonomy groupoid, and the map $\hat p_X:BG_X\to B\Gamma_n$ is the
one defined by \eqref{fF} in the standard way. According to
\cite{CM}, the composition of the map $\hat p_X^*$ with the
isomorphism \eqref{isom1} gives the map \eqref{CM19}. This, the
above construction and the diagram  imply that  \eqref{charmap}
coincides with  \eqref{CM19}.

Let us denote the map \eqref{newk} by $\tilde k$. We must prove that
$\tilde k$ coincides with $k$ given by \eqref{Bk}.

Let us first prove the equality
$$f^*_\F\circ \tilde k= f^*_\F\circ  k.$$

For the moment let us assume that $X$ is an arbitrary $D_n$-space.
Note that the domain of each chart from  the $\M_\infty$-atlas on
$S''(X)$ is a domain of the corresponding chart on $S''(\pt)$, and a
morphism of the charts on $S''(X)$ is a morphism of the
corresponding charts on $S''(\pt)$. This gives the homomorphism
$$(S''(p_X))^*:\vH_{\V_{S''}}^*(S''(\pt))\to \vH_{ \U_{S''}}^*(S''(X)).$$
From the definitions of the maps $j$ and $(S''(p_X))^*$ it follows
that we have the commutative diagram \be\label{diagr1}\xymatrix{
H^*(S''(\pt)) \ar@{->}^j[r]
\ar@{->}_{(S''(p_X))^*}[d]& \vH_{\V_{S''}}^*(S''(\pt))\ar@{->}^{(S''(p_X))^*}[d]\\
H^*(S''(X)) \ar@{->}^j[r] & \vH^*_{ \U_{S''}}(S''(X))}\ee

It is clear, that the isomorphism from Proposition~\ref{prop1} gives
the following commutative diagram: \be\label{diagr2} \xymatrix{
\vH_{\V_{S''}}^*(S''(\pt)) \ar@{<-}^\sim[r]
\ar@{->}_{(S''(p_X))^*}[d]& \vH_{\V}^*(\pt)\ar@{->}^{(p_X)^*}[d]\\
\vH^*_{ \U_{S''}}(S''(X)) \ar@{<-}[r]^{\,\,\sim} & \vH^*_{
\U}(X)}\ee Let now again $X=M/\F$. Let $\U_M$ be a cover of $M$
obtained from a maximal foliated atlas on $M$. Let $\U$ be the cover
of $X$ obtained from the full $\D_n$-atlas of $X$ defined by the
foliated atlas of $M$. The cover $\U_M$ gives the cover
$\U_{S''(\F)}$ of the manifold $S''(\F)$. The projection $p:M\to X$
provides the commutative diagrams
$$\xymatrix{S''(X) \ar@{->}[r]
\ar@{<-}_{S''(p)}[d]& X\ar@{<-}^{p}[d]\\
S''(\F) \ar@{->}[r] & M}$$ and $$ \xymatrix{\vH^*_{\U_{S''}}(S''(X))
\ar@{<-}[r]
\ar@{->}_{(S''(p))^*}[d]& \vH^*_\U(X)\ar@{->}^{p^*}[d]\\
\vH^*_{\U_{S''(\F)}} (S''(\F)) \ar@{<-}[r] & \vH^*_{\U_M}(M)}$$
Since both $M$ and $S''(\F)$ are smooth manifolds, the maps $j$
become isomorphisms and we get (recall that $S''(\F)$ and $M$ are
homotopy equivalent)
$$ \xymatrix{ \vH^*_{\U_{S''(\F)}}
(S''(\F)) \ar@{<-}[r]  \ar@{<-}^{\wr}[d]  & \vH^*_{\U_M}(M) \ar@{<-}^{\wr}[d]\\
H^*(S''(\F)) \ar@{<-}^\sim[r] & H^*(M) }
$$
Summarizing the above diagrams and using the isomorphisms
$H^*(\WOn)\cong H^*(S''(\pt))$ and \eqref{vHpt}, we obtain the
following diagram:
$$\xymatrix{ H^*(\WOn) \ar@{->}[r] \ar@{->}[d]&
\vH_{\V_{S''}}^*(S''(\pt))\ar@2{-}[r] \ar@{->}[d] &
H^*(B\Gamma_n)\ar@{->}[d]\\
H^*(S''(X)) \ar@{->}[r] \ar@{->}[ddr] & \vH^*_{ \U_{S''}}(S''(X))
\ar@2{-}[r] \ar@{->}[d]& \vH^*_{ \U}(X) \ar@{->}[d]\\
& \vH^*_{ \U_{S''(\F)}}(S''(\F)) \ar@2{-}[r] \ar@2{-}[d]& \vH^*_{
\U_M}(M) \ar@2{-}[d]\\
& H^*(S''(\F)) \ar@2{-}[r] & H^*(M) }$$

A part of the sequence of the maps of the diagram coincides with
\eqref{map3} and it gives the map $f^*_F \circ k$. The top row of
the maps give the map $\tilde k$. Since the maps \eqref{newk} and
\eqref{Bk} coincide, according to \cite{CM}, the right column of the
maps coincides with $f^*_F$. Thus, $f^*_\F\circ \tilde k=
f^*_\F\circ  k=k_\F.$

By Bott's construction \cite{Bott}, the map $k$ is uniquely defined
by the property $f^*_\F\circ  k=k_\F$ for all foliations $\F$. This
implies the equality $\tilde k=k$. The theorem is proved. \qed

\begin{theorem}\label{Thcompariz} Let $X$ be a $\D_n$-space with a full atlas $\Phi$
and the corresponding cover $\U$. Then there exists the following
commutative diagram
$$\xymatrix{
H^*(W_n,\On) \ar@{->}^k[r]
\ar@{->}_{\chi''}[d]& H^*(B\Gamma_n)\ar@{->}^{\check k_X}[d]\\
H^*(S''(X)) \ar@{->}^\alpha[r] & \vH^*_\U(X)}$$
\end{theorem}

{\bf Proof.} The proof of the theorem deliver the diagrams
\eqref{diagr1} and \eqref{diagr2}, and the isomorphisms
$H^*(W_n,\On)\cong H^*(S''(\pt))$ and \eqref{vHpt}.
 \qed

\begin{cor}\label{Corcompariz} Let $\F$ be a foliation of codimension $n$ on a smooth
manifold $M$. Let $T$ be a complete transversal for $T$ and $\U$ a
basis of $T$.  Then there exists the following commutative diagram:
$$\xymatrix{
H^*(\WOn) \ar@{->}^k[rr]  \ar@{->}_{\chi''}[d]
& & H^*(B\Gamma_n)\ar@{->}^{\check k_\F}[d]\\
H^*(S''(M/\F)) \ar@{->}^\alpha [rr]  \ar@{->}[rd]  & &
\vH^*_\U(M/\F) \ar@{->}[ld]
\\
& H^*(M) &}$$ The images of the generators of $H^*(\WOn)$ under the
map $\chi''$ are characteristic classes defined by Losik. The images
of the generators of $H^*(\WOn)$ under the map $\check k_\F\circ k$
are characteristic classes defined by Crainic and Moerdijk. These
classes are projected to the usual ones living in  $H^*(M)$.
\end{cor}

\begin{ex} Let $X$ be a  $\D_1$-space. The Godbillon-Vey class with values in $\vH^3(X)$ is obtained form the generator $[y_1c_1]\in H^3(W_1,{\rm O}(1))$, it is defined by the form of type $(2,1)$ given  by the formula \begin{equation}\label{formgv}gv(h_1,h_2)= \ln|h'_1|d\ln|h_2'\circ h_1|.\end{equation}
	For foliations of codimension one this formula is given in \cite{CM}. 
	
The Godbillon-Vey class with values in  $H^3(S''(X))$ may be described in the following way. Let $U\subset\Real$ be
the domain of a chart on $X$. A map $h$ from a neighborhood of $0$ in $\Real$ to $U$ defines the coordinates $z_i=h^{(i)}(0)$, $z=0,1,2,\dots$, on $S(U)$. Consider the following coordinates on $S''(U)$: $x_0=z_0,x_1=\ln|z_1|,x_k=z_k/z^k_1,$ $k=2,3,\dots$.
The Godbillon-Vey class is given by the form
$$dx_0\wedge dx_1\wedge dx_2$$ with respect to each such coordinate system.

The formula \eqref{formgv} defines a form of type $(2,1)$ on $S''(X)$ if the arrows
are changed to their extensions to $S''(X)$.
It is easy to check that
$$dx_0\wedge dx_1\wedge dx_2=gv+D(\omega^{(0,2)}+\omega^{(1,1)}+\omega^{(2,0)}),$$
where 
 $D$ is the total differential and
 $$\omega^{(0,2)}=x_1dx_2\wedge dx_1,\quad \omega^{(1,1)}(h)=x_2\ln|h'(h_0)|dx_0,\quad\omega^{(2,0)}(h_1,h_2)=\ln|h'_1|\ln|h_2'\circ h_1|.$$
This confirms the fact that the homomorphism $H^*(S''(X))\to\vH^*(X)$ respects the corresponding Godbillon-Vey classes.

Let now $\F$ be a Reeb
foliation on the three dimensional sphere. It is known that the
usual Godbillon-Vey class of $\F$ is trivial in $H^3(S^3)$. 
As it is explained in the Introduction, this class is also trivial in $\vH^3(S^3/\F)$. In
\cite{BGG} we show that for some specific choice of $\F$, the Godbillon-Vey class  is
non-trivial in $H^3(S(S^3/\F)/\O1)$. Hence, it allows to detect non-diffeomorphic foliations.
\end{ex}

\section{Characteristic classes of groups diffeomorphisms}

Let $G$ be an (abstract) group acting by diffeomorphisms on an $n$-dimensional manifold $M$. Let $H^*(G;\Omega(M))$ be the cohomology of the group $G$ with the coefficients in the forms on $M$ (which is defined in the same way as the \v{C}ech-de~Rham cohomology with all arrows being the global diffeomorphisms of $M$ from $G$).  The orbit space $M/G$ is a $\D_n$-space. It is obvious that there is a homomorphism
$$\vH^*(M/G)\to H^*(G;\Omega(M)).$$ From the above we obtain the sequence of maps
$$H^*(W_n,\On)\to H^*(S''(M/G))\to \vH^*(M/G)\to  H^*(G;\Omega(M)).$$
This gives the chracteristic map 
$$H^*(W_n,\On)\to  H^*(G;\Omega(M)).$$
 

\begin{ex} Let $G$ be a group of orientation preserving diffeomorphisms of the circle $S^1$. The formula \eqref{formgv}, for $h_1,h_2\in G$,  gives the Godbillon-Vey class with values in $H^3(G;\Omega(S^1)).$ 
	The integration of $gv(h_1,h_2)$ over $S^1$ gives the Bott-Thurston formula for the Godbillon-Vey class in $H^2(G;\Real)$, \cite{BottActions}. An important Theorem by Duminy and Sergiescu \cite{D-S} states that  if $G$ does not contain crossed elements, then this class is trivial. 
	
	On the other hand, let $\F$ be a Reeb foliation with non-trivial Godbillon-Vey class in $H^3(S^3/\F)$ constructed in \cite{BGG}.  Let $\xi$ be the product of two generators of the holonomy group of the compact leaf such that $\xi(x)<x$. Then $\xi$ is a diffeomorphism on its image of a neighborhood of $0\in\Real$ and  $0$ is the only fixed point of $\xi$. It is clear that $\xi$ may be extended to a diffeomorphism of $\Real$ satisfying $\xi(x+1)=\xi(x)+1$ and with the set $\mathbb{Z}$ of fixed points, i.e., $\xi$ defines a diffeomorphism of the circle with exactly one fixed point. This diffeomorphism generates a group $G$. It is clear that $G$ does not contain crossed elements. From the results of \cite{BGG} it follows that the Godbillon-Vey class of the obtained action is non-trivial in $H^3(S''(S^1/G))$. This shows that the Godbillon-Vey class in $H^3(S''(S^1/G))$ may be used to detect non-conjugate group actions.

\end{ex}

\section{Chern classes}

Till now, we considered only the cohomology $H^*(W_n,\On)$. It is
clear that in the same way we may consider the cohomology
$H^*(W_n,\GLn)$. Then the Proposition \ref{prop1}, Theorem
\ref{Thcompariz} and Corollary \ref{Corcompariz} hold true with
$S''(X)$ replaced by $S'(X)$. Moreover, there are the obvious maps
$H^*(W_n,\GLn)\to H^*(W_n,\On)$ and $S''(X)\to S'(X)$. For any
$\D_n$-space we obtain the diagram
$$\xymatrix{H^*(W_n,\GLn) \ar@{->}[r]
\ar@{->}_{\chi'}[d]& H^*(W_n,\On) \ar@{->}^k[r]
\ar@{->}_{\chi''}[d]& H^*(B\Gamma_n)\ar@{->}^{\check k_X}[d]\\
H^*(S'(X)) \ar@{->}[r] &H^*(S''(X)) \ar@{->}^\alpha[r] &
\vH^*_\U(X)}$$

\begin{cor}\label{corA} Let $X$ be a $\D_n$-space. Then the characteristic classes defined by the elements of the kernel of the homomorphism
$H^*(W_n,\GLn)\to H^*(W_n,\On)$ are zero in $\vH^*_{\U}(X)$. \end{cor}

\begin{ex}\label{ex2} Let $X$ be a $\D_1$-space.
The class $[c_1]\in H^2(W_1,{\rm GL}(1,\Real))$ defines the first
Chern class $C_1\in H^2(S'(X))$. If $x_0,x_1,x_2,\dots$ are the
standard coordinates on $S(U)$, then
$y_0=x_0,y_2=\frac{x_2}{x_0^2},y_3=\frac{x_3}{x_0^3},\dots$ are
coordinates on $S'(U)$. With respect to these coordinates, $C_1$ is
given by the form $$dy_2\wedge dy_0.$$ Let $X=M/\F$. The image of
the generator $[c_1]$ of $H^2(W_1,{\rm GL}(1,\Real))$ under the map
$H^2(W_1,{\rm GL}(1,\Real))\to H^2(W_1,{\rm O}(1))$ is trivial,
consequently the first Chern class is always trivial in the
\v{C}ech-de~Rham cohomology $\vH^2(M/\F)$ and in $H^2(M)$. This is
also directly proved in \cite{CM}. On the other hand, Losik
\cite{L15} (see also \cite{Baz-Gal})) proved that the first Chern class of the Reeb foliation
on the three-dimensional sphere is non-trivial in the cohomology
$H^2(S'(M/\F))$.\end{ex}

\begin{ex} Let $\F$ be a foliation of codimension one on a
three-dimension manifold. Suppose that the foliation $\F$ is defined
by a non-vanishing $1$-form $\omega$. Consider the complex
$$A^{m}=\Omega^{m-1}(M)\wedge \omega$$ with the differential being the usual exterior derivative. The cohomology of this complex are
denoted by $H^*_\F(M)$. Let $\eta$ be any 1-form such that
$d\omega=\eta\wedge \omega$. Then the Vey class is the class of the
form $d\eta$ in $H^2_\F(M)$. A Riemannian matric $g$ on $M$ defines
a map $\sigma:M\to S'(M/F)$ such that $\sigma^*(dy_2\wedge
dy_0)=d\eta$, where $\eta=-L_X\omega$ and $X$ is the vector field
orthogonal with respect to $g$ to the distribution tangent to $\F$
and such that $\omega(X)=1$. Unfortunately, the map $\sigma^*$ does
not induce a map from $H^*(S'(M/\F))$ to $H^*_\F(M)$. If $\F$ is the
Reeb foliation, then $d\eta=0$ \cite{Tamura}, i.e., the Vey class of
the Reeb foliation is trivial. 
\end{ex}

\section{Classes from the cohomology $H^*(W_n)$}

Let $X$ be a $\D_n$-space with a full atlas $\Phi$, and $\U$ the
corresponding cover. The proof of the following proposition is
similar to the proof of Proposition~\ref{prop1}.

\begin{prop} The projection $S(X)\to P(X)$ induces the isomorphism $$\vH^*_{\U_S}(S(X))\cong \vH^*_{\U_P}(P(X)).$$ \end{prop}

Let $X$ be a $\D_n$-space with a full atlas $\Phi$ that defines a
complete cover $\U$. Let $Q(X)$ be one of the spaces $S(X)$,
$S'(X)$, $S''(X)$, $P(X)$. Note that the cover $\U_Q$ is not
complete. Denote by $\bar\U_Q$ the cover obtaining as a base of the
topologies on $Q(U)$ for all $U\in\U$. The cover $\bar\U_Q$ is
complete. We need a complete cover in order to use the isomorphism
\eqref{lab*}. Let $U,V\in \U$. Let $\tilde U\subset Q(U)$ and
$\tilde V\subset Q(V)$ be open subsets. We will consider the
morphisms $g:\tilde U\to\tilde V$ that are restrictions of the
extensions $Q(h):Q(U)\to Q(V)$ of the morphisms $h:U\to V$ from
$\C_\Phi$. This gives a full atlas of $Q(X)$ and the corresponding
complete cover $\bar\U_Q$. The proof of the following proposition
will be given in Appendix.

\begin{prop}\label{propcomplcover} There is a natural isomorphism
$$\vH^*_{\U_Q}(Q(X))\cong \vH^*_{\bar\U_Q}(Q(X)).$$\end{prop}

Consider the frame bundle $P(\pt)$ for the space $\pt$. It is
obvious that the reduced groupoid $G_{P(\pt)}$ coincides with the
groupoid $\bar\Gamma_n$. The isomorphism \eqref{lab*} applied to
$P(\pt)$ gives \be\label{vHPpt} \vH^*_{\bar\V_P}(P(\pt))\cong
H^*(B\bar\Gamma_n).\ee

 Let $X$ be a $\D_n$-space. The
projection $p_X:X\to \pt$ induces the map
\be\label{newsecondary}H^*(W_n)\stackrel{\chi}{\longrightarrow}
H^*(S(X)) \stackrel{j}{\longrightarrow} \vH^*_{\U_S}(S(X))\cong
\vH^*_{\U_P}(P(X))\cong \vH^*_{\bar\U_P}(P(X)).\ee Taking $X=\pt$,
we get the map \be\label{newbark}\bar k:H^*(W_n)\to
H^*(B\bar\Gamma_n).\ee Applying to the projection $p_X:X\to \pt$ the
functor $P$, we get the maps $$P(p_X):P(X)\to P(\pt)$$ and
\be\label{HPXPpt} \bar H^*(B\bar\Gamma_n)\cong
\vH^*_{\bar\V_P}(P(\pt))\to \vH^*_{\bar\U_P}(P(X)).\ee

 We say that the
frame bundle $P(X)$ is trivializable if there is a section $s:X\to
P(X)$, i.e.,  for the domain $U$ of each chart form $\Phi$, there
exists a section $s_U:U\to P(U)$ such that for each morphism of
charts $m:U\to V$, the diagram
$$\xymatrix{P(U) \ar@{->}^{P(m)}[r]
\ar@{<-}_{s_U}[d]& P(V) \ar@{<-}^{s_V}[d]\\ U \ar@{->}_m[r] & V }$$
is commutative. It is clear that this notion does not depend of the
choice of a full atlas. It is obvious that if $X=M/\F$ is the leaf
space of a foliation $\F$ on a manifold $M$, then the bundle
$P(X)\to X$ is trivializable if and only if the bundle $P(\F)\to M$
is trivializable.

A section $s:X\to P(X)$ induces the map
$$\vH^*_{\bar\U_P}(P(X))\cong\vH^*_{\U_P}(P(X))\to \vH^*_{\U}(X).$$
Together with \eqref{newsecondary} and \eqref{HPXPpt} this delivers
the characteristic homomorphisms \be\label{WnvHX} H^*(W_n)\to
\vH^*_{\U}(X)\ee and \be\label{vbarfX} \check{\bar
f}_X^*:H^*(B\bar\Gamma_n)\to \vH^*_{\U}(X).\ee By the construction, the
map \eqref{WnvHX} coincides with the composition $\check{\bar
f}_X^*\circ \bar k$. This  defines the exotic characteristic
classes of a $D_n$-space $X$ with values in the \v{C}ech-de~Rham
cohomology of $X$. We obtain

\begin{theorem} Let $X$ be a $D_n$-space with a trivializable bundle
of frames, then there exist a characteristic map $\check{\bar
f}_X^*$ defined by \eqref{vbarfX} and a characteristic map
\eqref{WnvHX} that coincides with the composition $\check{\bar
f}_X^*\circ \bar k$.\end{theorem}

If $X=M/\F$ is the leaf space of a foliation $\F$ on a smooth
manifold $M$ with a trivializable bundle $P(\F)\to M$, then there is
the following commutative diagram:
$$\xymatrix{
H^*(W_n)\ar@{->}[r] & H^*(S(M/F)) \ar@{->}[r] & \vH^*_{\U}(M/F)\\
& H^*(M) \ar@{<-}[ur] \ar@{<-}[u] \ar@{<-}[ul] & }$$ Showing the
relation of the Losik  classes from $H^*(S(M/\F))$, the secondary
classes with values in the \v{C}ech-de~Rham cohomology of the leaf
space and the usual secondary characteristic classes of the
foliation $\F$.

The proof of the following theorem is similar to the proof of
Theorem \ref{Thcoincide}

\begin{theorem}\label{Thcoincide2}
The homomorphism \eqref{newbark} coincides with the homomorphism
\eqref{Bbk}.\end{theorem}

Completing this section note that the projection $S(\pt)\to
S''(\pt)$ together with the maps $j$ give the commutative diagram
$$\xymatrix{ H^*(S''(\pt)) \ar@{->}[r]
\ar@{->}[d]& H^*(S(\pt)) \ar@{->}[d]\\
\vH_{\V_{S''}}^*(S''(\pt)) \ar@{->}[r] & \vH^*_{ \V_{S}}(S(\pt))}$$
Applying the above described isomorphisms, we get the well-known
diagram
$$\xymatrix{ H^*(\WOn) \ar@{->}[r]
\ar@{->}[d]& H^*(\Wn) \ar@{->}[d]\\
H^*(B\Gamma_n) \ar@{->}[r] & H^*(B\bar\Gamma_n)}$$ Note that the map
$H^*(B\Gamma_n)\to H^*(B\bar\Gamma_n)$ is in fact induced by the
projection $P(\pt)\to \pt$.

\section{The full picture}

Finally, considering the projections $S(X)\to S''(X)\to S'(X)$, we obtain the following theorem.

\begin{theorem} Let $X$ be a $\D_n$-space with a section $s:X\to P(X)$. Then there exists the following commutative diagram

\newpage

\vskip1.4cm

$$\xymatrix{
H^*(B\bar\Gamma_n)\ar@{<-}[d] \ar@{<-}[r] \ar@/^4pc/[rrrrd]& H^*(W_n)\ar@{->}^\chi[r] \ar@{<-}[d]            & H^*(S(X))\ar@{<-}[d]   \ar@{->}[rrd] \ar@{<-}[d]       &    &  \\
H^*(B\Gamma_n) \ar@{<-}[r]  \ar@/_6.5pc/[rrrr] & H^*(W_n,\On) \ar@{->}^{\chi''}[r] \ar@{<-}[d]   & H^*(S''(X)) \ar@{<-}[d]      \ar@{->} [rr] \ar@{<-}[d]  &    &  \vH^*_{\U}(X)  \\
& H^*(W_n,\GLn) \ar@{->}^{\chi'}[r]               & H^*(S'(X))
\ar@{->} [urr]  & }$$

\vskip1.4cm

Here the images of the maps $\chi$, $\chi'$ and $\chi''$ give the
characteristic classes defined by Losik, and the images of the maps
from $H^*(W_n,\On)$ and $H^*(W_n)$ to $\vH^*_{\U}(X)$ give the
 charcteristic classes with values in the
\v{C}ech-de~Rham cohomology of $X$.

\end{theorem}

Note that the map $H^*(W_n,\GLn)\to H^*(W_n)$ is zero. This implies

\begin{cor} Let $X$ be a $\D_n$-space with a section $s:X\to P(X)$. Then the Chern classes, i.e., the characteristic classes defined by the cohomology $H^*(W_n,\GLn)$,
 are zero in~$\vH^*_{\U}(X)$. \end{cor}

We also get.

\begin{cor} Let $X$ be a $\D_n$-space with a section $s:X\to P(X)$. Then the characteristic classes defined by the elements of the kernel of the homomorphism
$H^*(W_n,\On)\to H^*(W_n)$ are zero in $\vH^*_{\U}(X)$. \end{cor}

Let $n=1$. The homomorphism $H^*(W_1,{\rm O}(1))\to H^*(W_1)$ is an isomorphism.
The last two corollaries do not give new information comparing to Corollory~\ref{corA} and Example~\ref{ex2}.

Let $n=2$. The element $c_2$ is trivial in $H^*(W_2)$, and it is
non-trivial in both $H^*(W_2,\GLn)$ and $H^*(W_2,\On)$.
Consequently, if  $X$ is a $\D_2$-space with trivial $P(X)$, then
the image of $[c_2]$ is trivial in $\vH^*_{\U}(X)$. It would be
interesting to construct a foliation of codimension two with
trivializable $P(M/\F)$ and such that $[c_2]$ is non-trivial in
$H^4(S''(M/\F))$ (and consequently non-trivial in $H^4(S'(M/\F))$).

\appendix
\section{Proof of Proposition \ref{propcomplcover}}

Recall that the cochains form $\check C^{k,q}_{\bar \U_Q}(Q(X))$ map
 the strings of morphisms
$$V_0\stackrel{g_1}{\longrightarrow}\cdots
\stackrel{g_k}{\longrightarrow} V_k$$ to $\Omega^q(V_0)$, here each
$V_i$ is an open subset of some $Q(U_i)$, $U_i\in\U$,
$g_i:V_{i-1}\to V_i$ is the restriction to $V_{i-1}$ of the
extension $Q(h_i):Q(U_{i-1})\to Q(U_i)$ of a morphism of charts
$h_i:U_{i-1}\to U_i$. It is clear that we may assume that
$p(V_i)=U_i$, where $p:Q(U_i)\to U_i$ is the projection. The
cochains form $\check C^{k,q}_{\U_Q}(Q(X))$ map
 the strings of morphisms
$$Q(U_0)\stackrel{Q(h_1)}{\longrightarrow}\cdots
\stackrel{Q(h_k)}{\longrightarrow} Q(U_k)$$ to $\Omega^q(Q(U_0))$.
We define the morphisms of the complexes
  $$\mu:\check C^{k,q}_{
\U_Q}(Q(X))\to \check C^{k,q}_{\bar \U_Q}(Q(X)),$$ $$\lambda:\check
C^{k,q}_{\bar \U_Q}(Q(X))\to \check C^{k,q}_{ \U_Q}(Q(X))$$ by
setting
$$\mu(\varphi)\left(V_0\stackrel{g_1}{\longrightarrow}\cdots
\stackrel{g_k}{\longrightarrow}
V_k\right)=\varphi\left(Q(U_0)\stackrel{Q(h_1)}{\longrightarrow}\cdots
\stackrel{Q(h_k)}{\longrightarrow} Q(U_k)\right)\Big|_{V_0},\quad
\varphi\in \check C^{k,q}_{ \U_Q}(Q(X)),$$
$$\lambda(c)\left(Q(U_0)\stackrel{Q(h_1)}{\longrightarrow}\cdots
\stackrel{Q(h_k)}{\longrightarrow}
Q(U_k)\right)=c\left(Q(U_0)\stackrel{Q(h_1)}{\longrightarrow}\cdots
\stackrel{Q(h_k)}{\longrightarrow} Q(U_k)\right),\quad c\in \check
C^{k,q}_{\bar \U_Q}(Q(X)).$$
 It holds
$$\lambda\circ\mu=\id,$$ i.e. $\mu$ induces a monomorphism in
cohomology. We are going to construct a chain homotopy between the
maps $\mu\circ\lambda$ and $\id$. We define the map
$$F:\check
C^{k,q}_{\bar \U_Q}(Q(X))\to C^{k-1,q}_{\bar \U_Q}(Q(X))$$ in the
following way:
\begin{multline*}F(\varphi)\left(V_0\stackrel{g_1}{\longrightarrow}\cdots
\stackrel{g_{k-1}}{\longrightarrow}
V_{k-1}\right)\\=\sum_{s=0}^{k-1}(-1)^s
\varphi\left(V_0\stackrel{g_1}{\longrightarrow}\cdots
\stackrel{g_s}{\longrightarrow} V_{s}
\stackrel{i_s}{\longrightarrow}
Q(U_s)\stackrel{Q(h_{s+1})}{\longrightarrow}\cdots
\stackrel{Q(h_{k-1})}{\longrightarrow} Q(U_{k-1})
\right),\end{multline*} where $i_s:V_{s}\to Q(U_s)$ is the
inclusion. It is clear that $F$ commutes with the differential $d$.
Next, for the total differential $D$ and $\varphi\in C^{k,q}_{\bar
\U_Q}(Q(X))$ it holds
\begin{multline*}D(F(\varphi))+F(D(\varphi))=(\delta+(-1)^{k-1}d)(F(\varphi))+F((\delta+(-1)^{k}d)\varphi)\\=\delta(F(\varphi))+F(\delta\varphi).\end{multline*}
Hence in order to show that $F$ is a cochain homotopy it is
sufficient to prove the equality
\be\label{cheinhom}\mu\circ\lambda-\id=\delta\circ
F+F\circ\delta.\ee This equality may be checked directly. To
illustrate this let us for simplicity suppose that $\varphi\in
C^{2,q}_{\bar \U_Q}(Q(X)).$ It holds
$$(\mu\circ\lambda-\id)(\varphi)(g_1,g_2)=\varphi(Q(h_1),Q(h_2))|_{V_0}-\varphi(g_1,g_2).$$
Next, \begin{multline*} F(\delta
\varphi)(g_1,g_2)=(\delta\varphi)(i_0,Q(h_1),Q(h_2))-(\delta\varphi)(g_1,i_1,Q(h_2))+(\delta\varphi)(g_1,g_2,i_2)\\=
\varphi(Q(h_1),Q(h_2))|_{V_0}-\varphi(Q(h_1)\circ
i_0,Q(h_2))+\varphi(i_0,Q(h_2)\circ
Q(h_1))-\varphi(i_0,Q(h_1))\\
-(g_1^*\varphi(i_1,Q(h_2))-\varphi(i_1\circ
g_1,Q(h_2))+\varphi(g_1,Q(h_2)\circ i_1)-\varphi(g_1,i_1))\\
+g_1^*\varphi(g_2,i_2)-\varphi(g_2\circ
g_1,i_2)+\varphi(g_1,i_2\circ g_2)-\varphi(g_1,g_2);
\end{multline*}
\begin{multline*}
\delta(F(\varphi))(g_1,g_2)=g_1^*F(\varphi)(g_2)-F(\varphi)(g_2\circ
g_1)+F(\varphi)(g_1)\\
=g_1^*\varphi(i_1,Q(h_2))-g_1^*\varphi(g_2,i_2)-\varphi(i_0,Q(h_2\circ
h_1))+\varphi(g_2\circ
g_1,i_2)+\varphi(i_0,Q(h_1))-\varphi(g_1,i_1).
\end{multline*}
Noting that $Q(h_{s+1})\circ i_s=i_{s+1}\circ g_{s+1}$, we see that
\eqref{cheinhom} holds true for the $\varphi$ under consideration.
The proposition is proved. \qed

 \vskip1cm

{\bf Acknowledgements.} The author is thankful to Steven Hurder for
useful email communications, and to Yaroslav Bazaikin and Alexei
Kotov for fruitful discussions. The work was supported by grant no.
18-00496S of the Czech Science Foundation.

\end{document}